\documentclass[a4paper]{article}

\usepackage{natbib}
\usepackage{graphicx}
\usepackage{amssymb}
\usepackage{url}

\newcommand{\latin}[1]{\textit{#1}}  
\newcommand{\firstuse}[1]{\emph{#1}} 
\newcommand{\email}[1]{\texttt{#1}} 
\newcommand{\ie}{i.e.\ }

\newcommand{\unix}[1]{\texttt{#1}} 
\newcommand{\card}[1]{\vert#1\vert} 
\newcommand{\eclipse}{ECL$^i$PS$^e$\ }

\providecommand{\keywords}[1]
{
  \small
  \textbf{\textit{Keywords---}} #1
}

\begin{document}

\title{Selecting and scheduling an optimal subset of road network upgrades}

\author{A.D.~Stivala\thanks{NICTA Victoria Research Laboratories, Australia} \thanks{Department of Computer Science and Software Engineering, The University of Melbourne, Victoria 3010, Australia}
 \and
P.J.~Stuckey\footnotemark[1] \footnotemark[2]
 \and
M.G.~Wallace\footnotemark[1] \thanks{Faculty of Information Technology, Monash University, Caulfield East, Victoria 3145, Australia. Email: \email{mark.wallace@monash.edu}}
}

\date{November 15, 2011}

\maketitle

\begin{abstract}
We consider the problem of choosing a subset of proposed
road network upgrades to implement within a fixed budget in order
to optimize the benefit in terms of vehicle hours travelled (VHT),
and show how to render the solution of this problem more tractable by
reducing the number of traffic assignment problems that must
be solved. This technique is extended to develop a heuristic algorithm
for finding a schedule of road upgrades over a planning period that maximizes the net present value of the resulting VHT reductions.
\end{abstract}



\section{Introduction}
\label{sec:intro}

Consider the problem of deciding, amongst a set of proposed upgrades
to an urban road network, which subset to implement given a fixed
budget that constrains the selection. In order to make this selection
in a rational way, we must associate some benefit with each possible subset.
One of the most obvious such benefits is simply the reduction in total travel
time across the system, which we can measure in \firstuse{vehicle hours
travelled (VHT)} per day, due to the upgrade. 

In order to compute the change in VHT for a particular upgrade, we
must solve the \firstuse{traffic assignment problem (TAP)}. This is an
important problem in transportation research and in the practical
application of modelling to transportation planning. This problem is
to predict the route choice of users of a road network given sets of
origins and destinations for trips. Once these routes have been
computed, the level of congestion on each road is determined, and
proposed changes to a road network can be compared by their effects on
the final congested network.

The importance of the traffic assignment problem has led to it being
an active area of research for many years, and a comprehensive
overview can be found in \citet{patriksson94}. Although there are many
approaches and variations to traffic modelling, we restrict ourselves
to the ``traditional'' \firstuse{static assignment} at \firstuse{user
  equilibrium (UE)}. In this model, we make the assumption that each
traveller chooses the route that minimizes his or her own travel time,
given the current traffic flows on the road network.  These travel
times are a function of the traffic flows resulting from the routes
chosen by all the travellers, and the user equilibrium is achieved
when no traveller can improve his or her travel time by unilaterally
choosing a different route.  
It is important to note that adding capacity (such as a new road) to a road
network does not necessarily reduce the total VHT, since the UE, 
which is optimal for each traveller, is not necessarily the optimum
for the system as a whole. Such a situation, when adding capacity
increases the system-wide travel time, is known as Braess's paradox.

This UE model has been formulated as a
convex quadratic optimization problem \citep{beckmann56,dafermos69}, which can be
solved efficiently by the Frank-Wolfe algorithm
\citep{frank56,leblanc75}. Although more sophisticated and faster converging algorithms
have been demonstrated, such as path-based assignment
\citep{jayakrishnan94,dial06}, origin-based assignment
\citep{bargera02}, and LUCE \citep{gentile09}, the simplicity of the
Frank-Wolfe algorithm, its small memory requirements, and its capacity
for simple parallelization mean that it is still competitive on
efficiency in practice, and widely used, for example in commercial
products such as PTV A.G.'s VISUM \citep{ptvag}
(which now also implements the LUCE algorithm), Citilabs's Cube
Voyager \citep{citilabs}, and Caliper's
TransCAD \citep{caliper,slavin06}.

In order to solve our constrained optimization problem (maximize
reduction in
VHT by selecting a subset of possible upgrades within a fixed budget),
we must at least compute the reduction in VHT due to each individual
upgrade in isolation. However, this might not be enough, since the
change in VHT due to two or more upgrades made to the network
need not necessarily simply be the sum of the changes individually.
To solve the problem exactly, we would in principle have to compute
the VHT for every member of the power set of the set of possible 
upgrades, an exponentially large number of traffic assignment problems
that is clearly infeasible in practice, especially given that each
TAP takes several hours to solve to acceptable accuracy for a realistic
city road network. By considering only pairwise interactions, on the assumption
that interactions of a higher order can be ignored, this can be reduced
to a quadratic number of TAPs to solve, which for a large number of possible
upgrades is still most likely impractical. In this paper we show
how to approximate an optimal solution within a practical time
by solving a model which incorporates pairwise interactions between
upgrades only when two individual road upgrades are \emph{estimated}
to have a pairwise VHT change value that is significantly different
from the sum of their individual contributions. That is, we need only
solve the TAP for each individual upgrade, and for a small number
of pairs of upgrades. Further, we present a model, and a heuristic algorithm
to approximate an optimal solution, for the problem of scheduling the upgrades
in a planning period in order to maximize the net present value of the VHT
reductions.

Recently, \citet{wu11} presented a formulation and solution of the
``project selection problem'' in the context of build-operate-transfer
developments of road networks. Their formulation involves finding 
projects that both have a societal benefit (VHT reduction) and 
are also attractive to private developers bidding for the development,
a function of the tolls they are able to charge. Therefore, their formulation
involves finding the optimal capacity increases and tolls for projects
to be contracted to private firms under the build-operate-transfer scheme.
In contrast, the project selection problem we solve here operates under
the simpler scheme of selecting from a set of candidate highway upgrades
those with the maximal benefit within a fixed budget. Our focus is on reducing
the number of TAPs to solve without unduly reducing the accuracy of the
resulting VHT reduction estimates. In addition, we extend the problem
to project scheduling, where we show a method to select the projects
to build in each year of a planning period, given a budget for each
year, in order to maximize the net present value of the societal 
benefit (reductions in VHT).

\section{Methods}

In order to find the optimal subset of a set of road upgrades within a
fixed budget, we need to first solve the TAP on the existing network
to obtain a VHT value for the current (before upgrade(s)) situation.
From then on we consider changes (hopefully reductions) in the VHT
from this base value for each upgrade. Solving the TAP requires the 
usual inputs to this process: the road network, with parameters
specifying the relationship between volume (flow) and travel time (cost)
on each link, and the \firstuse{Origin-Destination (O-D) demand matrix}
specifying trips on the network. The output of this process is then
the volume and latency (cost) on each link at UE, and global statistics such as VHT,
which is the sum over all links of the product of the volume
and cost on the link.

The outline of the procedure to find the optimal upgrade subset is then:

\begin{enumerate}
\item For each individual upgrade, compute VHT by solving the traffic assignment problem on the network modified by that upgrade.
\item Use the pairwise interaction model (Section~\ref{sec:pairwise_predictor})
  to estimate which pairs of upgrades will have a significant
  interaction, meaning that their combined effect is significantly
  different from simply additive.
\item Run traffic assignment on each such pair of upgrades to compute
  the VHT for these pairs of upgrades.
\item Solve the upgrade subset model (Section~\ref{sec:upgrade_subset_model}) with
  the individual and required pairwise upgrade VHT change values to find the
  optimum subset of upgrades within the specified budget.
\end{enumerate}

\subsection{Solving the traffic assignment problem}
\label{sec:tap}

As mentioned in Section~\ref{sec:intro}, the traffic assignment problem
is a well-studied problem.  Here we give a brief summary of the 
user equilibrium in the static traffic assignment problem,
and its solution by the Frank-Wolfe method. Details
can be found in, for example, \citet{sheffi85,patriksson94,florian95}.

We are given the road network as a graph $G = (\mathcal{N}, \mathcal{A})$
where $\mathcal{N}$ is the set of nodes and $\mathcal{A}$ is the set
of arcs (links representing road segments). Each arc $a \in \mathcal{A}$
has associated with it a capacity (the maximum number of vehicles
that can travel on it per unit time) $Q^{\max}_a$ and a free travel time
$t_a$, the time required to traverse it at the maximum allowed speed.
It is assumed that this travel time depends only on the current 
flow on the arc $f_a$ (the number of vehicles using that link), according
to a latency function which is continuous and non-decreasing (higher
flow means longer latency). The BPR (Bureau of Public Roads) function
is usually used for this purpose:
\begin{equation}
  \label{eqn:bpr}
  l_a(f_a) = t_a\left(1+\alpha\left(\frac{f_a}{Q^{\max}_a}\right)^\beta\right)
\end{equation}
where $\alpha$ and $\beta$ are parameters, for example $\alpha = 0.15$
and $\beta = 4.0$.

In addition, we are given a set of origin nodes $\mathcal{O} \subseteq \mathcal{N}$ and destination nodes $\mathcal{D} \subseteq \mathcal{N}$ and 
an O-D demand matrix $Q_{\card{\mathcal{O}}\times\card{\mathcal{D}}}$, 
such that $(q_{rs})$ gives the number of vehicles from origin $r \in \mathcal{O}$
to destination $s \in \mathcal{D}$ in a given time period. These origin
and destination nodes are typically the centroids of zones used in 
determining trip rates.

The traffic assignment which we are to find is the current state of
the traffic network given by the pair $(f,l)$ where $f \in
\mathbb{R}^{\card{\mathcal{A}}}$ is the flow (the volume of cars on
each link) and $l \in \mathbb{R}^{\card{\mathcal{A}}}$ is the travel
time (the latency on each link), where $l$ is a function of $f$
(Equation~\ref{eqn:bpr}).  We denote the flow on path $k$ connecting an O-D
pair $(r,s) \in \mathcal{O} \times \mathcal{D}$ by $h_{rs}^k$,
where $k$ is in $\mathbb{P}_{rs}$, the set of paths between origin
$r$ and destination $s$.

The user equilibrium is then obtained by solving the convex optimization problem:
\begin{equation}
  \label{eqn:tap_ue_obj}
  \min \sum_{a \in \mathcal{A}} \int_{0}^{f_a} l_a(x) dx
\end{equation}
subject to
\begin{eqnarray}
  \label{eqn:tap_ue_constraints}
  \sum_{k \in \mathbb{P}_{rs}} h_{rs}^k = q_{rs} && \forall (r,s) \in \mathcal{O}\times\mathcal{D}\\
  h_{rs}^k \geq 0 && \forall k \in \mathbb{P}_{rs}, \forall (r,s) \in \mathcal{O} \times \mathcal{D} \\
  f_a = \sum_{(r,s) \in \mathcal{O}\times\mathcal{D}} \sum_{k \in \mathbb{P}_{rs}} h_{rs}^k \delta_{a,k}^{rs} && \forall a \in \mathbb{A}\\
\end{eqnarray}
where $\delta_{a,k}^{rs}$ is the indicator variable
\begin{equation}
  \label{eqn:tap_ue_delta}
  \delta_{a,k}^{rs} = \left\{
    \begin{array}{lr}
      1, & \mbox{if link } a \mbox{ is on path } k \mbox{ between O-D pair } (r,s) \\
    0, & \mbox{otherwise} \\
    \end{array}
    \right.
\end{equation}

As mentioned in Section~\ref{sec:intro}, this quadratic program is often
solved with the Frank-Wolfe algorithm \citep{frank56} (also known as
the \firstuse{convex combinations} method).  Applied to the UE 
problem (Equations~\ref{eqn:tap_ue_obj}--\ref{eqn:tap_ue_delta}), this
amounts to iteratively solving a set of shortest-path problems
(the \firstuse{``all-or-nothing assignment''} of traffic flows) to find search
directions followed by a linear program to find search step size, terminating
when a convergence criterion is satisfied.

Many termination criteria are possible, but probably the most commonly
used is the \firstuse{relative gap} \citep{rose88,slavin06}.
This is the ratio of the 
difference between the cost of the current solution
and the cost of the ``all-or-nothing assignment''
solution (\ie the shortest paths that were used to find the search direction)
to the cost of the current solution:
\begin{equation}
\label{eqn:relgap}
\frac{\sum f \cdot l(f) - \sum f_{\mbox{aon}} \cdot l(f)}{\sum f \cdot l(f)}
\end{equation}
where $f_{\mbox{aon}}$ is the ``all-or-nothing'' solution.
Relative gaps on the order of $10^{-5}$ (or smaller) are usually considered 
to indicate a very good solution.

Since the shortest path computation is responsible for more than 90\%
of the CPU time in practical applications \citep[p. 100]{patriksson94},
an algorithm that solves this problem for road networks efficiently in practice
is important. Although Dijkstra's algorithm \citep{dijkstra59}
with an appropriate priority queue implementation such as a
Fibonacci heap \citep{fredman87} is among the 
algorithms with asymptotically smallest time complexity for the single-source
shortest path algorithm for general graphs with non-negative edge
weights, in practice on large-scale
road networks, \firstuse{label-correcting} algorithms run faster
\citep{zhan98,klunder06}. In contrast to the \firstuse{label-setting} algorithms,
such as Dijkstra's, which determine exactly the shortest path for one
node at each iteration, the label-correcting algorithms do not determine
shortest paths for any node until the end of the algorithm. That is,
the so-called pivot node (the node whose shortest path is being updated)
can be entered into the processing queue multiple times. Each time
it is dequeued its label value is an estimate (actually an upper bound) 
on the cost of the shortest path to that node from the source,
and only at termination (when the queue is empty because no smaller
labels are found) are the labels exact values for the shortest path costs.

The simplest label-correcting algorithm is the Bellman-Ford algorithm
\citep{bellman57,ford65}, however in practice on large-scale road
networks, better performance can be achieved with other algorithms
\citep{klunder06} such as the d'Esopo-Pape algorithm \citep{pape74},
and the \firstuse{Small Label First
  (SLF)} \citep{bertsekas93} algorithm, especially with the
 \firstuse{Large Label Last (LLL)}
\citep{bertsekas96} modification.
The d'Esopo-Pape algorithm always removes nodes
from the top of the queue, and enqueues nodes at the bottom of the
queue if they have never been in the queue before, otherwise at the
top of the queue.  The SLF algorithm is, instead of enqueuing
according to whether previously in queue or not, whenever a node is
enqueued, its label  (current distance) is compared with that at the top of the queue, and
if its label is smaller than that at the top, it is put at the top,
otherwise at the bottom of the queue.  In the LLL modification of the
algorithm, at each iteration, when the node at the top of the queue
has a larger label than the average label of nodes
in the queue, this (top) node is not removed as the current node but
instead put at bottom of queue (and we keep doing this until a node
with smaller than average label is found). 
\citet{klunder06} show that the d'Esopo-Pape and SLF algorithms
with the LLL modification
are faster in practice on large-scale road networks than
Dijkstra's algorithm  and the Bellman-Ford
algorithm.

\subsubsection{Implementation}

We implemented a Frank-Wolfe solution of the above traffic assignment
problem in the C programming language, using the 
d'Esopo-Pape algorithm with the LLL modification (we found that 
on our test network (see Section~\ref{sec:results}) the
d'Esopo-Pape algorithm with the LLL modification was faster
than the SLF algorithm with the LLL modification) to solve
the shortest path problems. Since the Frank-Wolfe algorithm involves
solving the shortest path problem for each O-D node pair at each
iteration, and each of these computations is independent, we can
execute these computations in parallel, and our implementation
does so using POSIX threads (\unix{pthreads}) to take
advantage of multicore processors. The resulting updates of the
link volumes vector are also executed in parallel, in a lock-free
manner using processor atomic instructions to ensure correctness.

We require the solution of the TAP not only for the original
road network, but also for the road network modified by the set of
potential upgrades, some pairs of potential upgrades, and, for 
evaluation purposes in our test cases, all pairs of upgrades and
all subsets of the upgrade set. Hence our implementation takes as
an additional input a file containing the description of the upgrade
set, with options to run with each modification, pairwise modifications,
or all subsets of modifications. Since these different modified 
road networks are processed independently, they can be run in parallel
and the Message Passing Interface (MPI) is used to distribute the
computations across different processors and nodes where available.

Source code for our implementation is available from 
\url{https://github.com/stivalaa/traffic_assignment}.

\subsection{A model to determine the optimal subset of upgrades}
\label{sec:upgrade_subset_model}

If, for $N$ upgrades, we know the VHT change for each individual upgrade, $v_i$,
and the VHT change for each pair of upgrades $\Delta\mbox{VHT}_{ij}$,
giving us the difference
between the pairwise change and the sum of individual changes,
$d_{ij} =  \Delta\mbox{VHT}_{ij} - (v_i + v_j)$
then the optimal subset can be found by the following constrained
optimization problem:
\begin{equation}
\label{eqn:objective}
\max \left( 
\sum_{1 \leq i < j \leq N} y_i y_j d_{ij} m + \sum_{1 \leq i \leq N} y_i (  v_i m  - c_i)
\right)
\end{equation}
subject to
\begin{equation}
  \label{eqn:budget_constraint}
  \sum_{1 \leq i \leq N} y_i c_i \leq B
\end{equation}
where $c_i$ is the cost of upgrade $i$ and $B$ is the total budget
and $y_i \in \{0,1\}$ is an indicator variable such that $y_i = 1$
if upgrade $i$ is in the selected subset.

In order to convert a VHT reduction $v$ to a monetary amount,
the factor $m$ is required that gives a dollar value to each vehicle hour
travelled. This may be some estimate of the average value of each traveller's
time, and will likely be some value on the order of \$10, for example.
In addition, since our traffic models operate
on daily VHT values, we must multiply by 365 to convert to annual
values, so, combining these considerations we might use for example $m = 3650$. 
We use this value of $m$ to compute the results presented in Section~\ref{sec:results}.
In a more detailed analysis we might have separate models for AM and PM
peaks on weekdays and for weekends, and hence the appropriate different
factors will be used to convert these to annual monetary benefits.

The objective function (Equation~\ref{eqn:objective})
is quadratic because as well as the benefit
for each individual upgrade we account for the possibly different 
benefit for each pair of upgrades in the case both in the pair are in
the subset considered.

Note the VHT change $v_i$ for a change $i$ is $\mbox{VHT}_0 - \mbox{VHT}_i$
where $\mbox{VHT}_0$ is the baseline VHT value, that is, the VHT for 
the existing road network with no upgrades. Therefore a positive value
of $v_i$ denotes an improvement (\ie reduction) in the VHT value
relative to the baseline due to change $i$.

This model can be easily expressed in the MiniZinc
\citep{nethercote07} modelling language and solved with the \eclipse
interval arithmetic constraint solver \citep{schimpf10,eclipse}.

If the VHT changes for all pairwise upgrades are computed, then all
values of $d_{ij}$ are known and this model finds the optimal subset
of upgrades exactly, up to pairwise interactions (note we have already
excluded the effect of interactions of a higher order).  However, this
is still a quadratic number of traffic assignment problems to
solve. On the assumption that only certain pairwise interactions need
to be considered, and that we can make some reasonable prediction of
which these are after solving the TAP for only each individual change,
we can simply set $d_{ij} = 0$ (\ie assume that $\Delta\mbox{VHT}_{ij}
= v_i + v_j$) for two changes $i$ and $j$ that are predicted to have
no significant pairwise interaction.

\subsection{Predicting upgrades with significant interactions}
\label{sec:pairwise_predictor}

When we consider a set of potential upgrades to a road
network, an intuitive assumption might be that upgrades in distant
parts of the network, particularly on ``local'' roads, will not affect
each other. Conversely, upgrades to the same road, or to different
roads that are very close together or intersect, or two alternatives
consisting of upgrading one route or adding an alternative route,
would be assumed to have a significant interaction.

For instance, consider the set of upgrades to the Chicago Regional
road network in Table~\ref{tab:chicago_mods}, 
illustrated in Figure~\ref{fig:chicago_mods}. With the exception
of the two extensions of Joe Orr Road (projects 07-94-0027 and 07-96-0013),
these projects are all in quite separate parts of the network and might
be assumed to be independent. We would certainly not make this assumption
about the two extensions to the same road however, and would want to
solve the TAP for both these changes together rather than just summing
the VHT reductions resulting from each independently.

Similarly with the list of hypothetical upgrades to the Berlin Center
road network listed in Table~\ref{tab:berlin_mods} and shown in
Figure~\ref{fig:berlin_mods_overview}, we would assume that 
the upgrades ber03, ber04, ber05, ber06, and ber06a might have significant
pairwise interactions, as well as the pair ber10 and ber10a, but that
all other upgrades could be assumed to be independent.

The results detailed in the next section show that these assumptions
are reasonable, and that accurate approximations can be made by including
only these pairwise interactions that are \latin{a priori} considered
significant. However, where there are a large number of proposed upgrades,
making it impractical to make these decisions ``by eye'', or an automatic
method is required for some other reason, at least two alternative
automatic techniques are possible.  The first is to simply choose some
threshold distance and consider all pairs of upgrades that are less
than this (Euclidean) distance apart to be not independent 
(\ie require a pairwise TAP computation).

For example, in the Chicago Regional road network and the upgrades
listed in Table~\ref{tab:chicago_mods}, the upgrade pair
07-94-0027 and 07-96-0013 has the smallest Euclidean distance
between the (centroids of) the upgraded links, an order of magnitude
smaller than the next smallest. Similarly, in the Berlin Center road network
and the upgrades in Table~\ref{tab:berlin_mods}, the 11 pairs of upgrades
that we have considered \latin{a priori} to have possibly significant
pairwise interaction coincide exactly with the smallest 11 in the sorted
list of Euclidean distances between (centroids of) changed links in
the upgrades.

The second technique is to 
use a clustering algorithm (such as $k$-means clustering) to cluster
the set of upgrades on their Euclidean $(x,y)$ co-ordinates and consider
all pairs within each cluster as potentially significantly interacting
(requiring pairwise TAP computations). In the latter case an algorithm
or parameters to the algorithm should be set so that there is a tendency
for the algorithm to find a large number of singleton clusters
in order to avoid unnecessarily clustering upgrades together that 
do not have significant pairwise interaction.

Other attributes such as network distance
(counting links rather than Euclidean distance), sums of flow changes on
links,  overlaps in the sets
of links with flow changes, and the number of shortest paths (and the flow
on them) in common
between two upgrades can be used as input to a clustering or classification
algorithm to determine pairs of upgrades that should be considered
as possibly having significant interaction. However we have found
no advantage in using these over the simpler distance-based threshold
or clustering techniques.

\section{Results}
\label{sec:results}

\begin{figure}
  \includegraphics[width=0.9\textwidth]{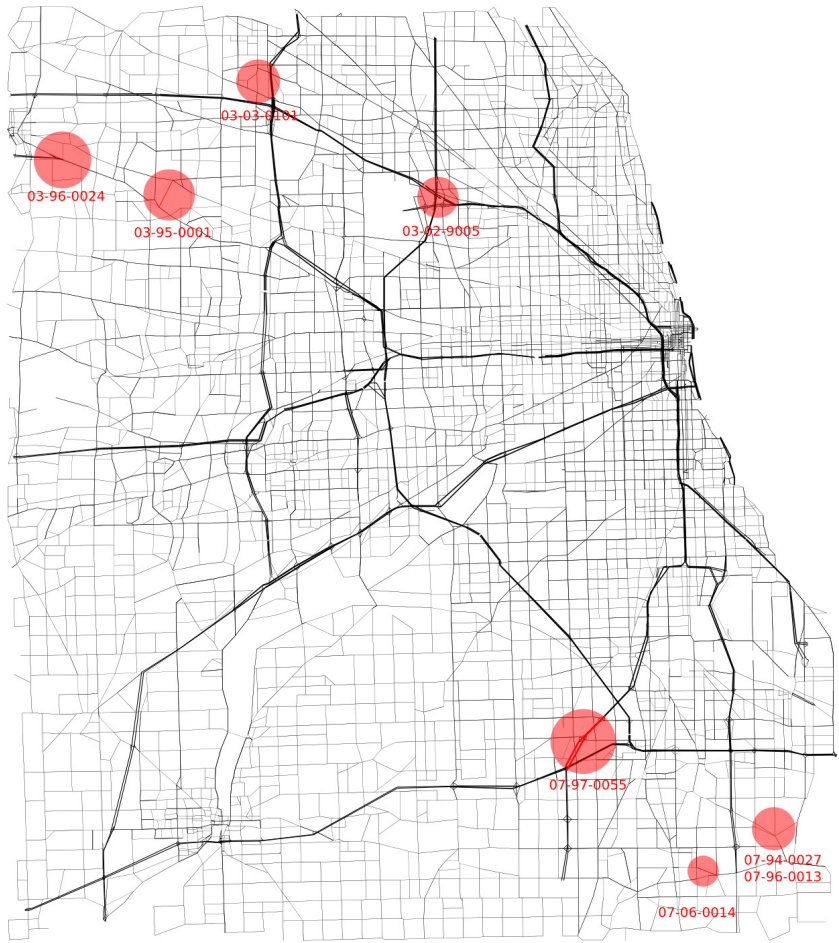}
  \caption{Potential upgrades to the Chicago Regional road network.}
  \label{fig:chicago_mods}
\end{figure}

We tested our methods on two of the realistic sized road networks from
Dr Hillel Bar-Gera's Transportation Test Problems website \citep{tntp}.
The first is the Chicago Regional road network.
This network,
provided by the Chicago Area Transportation Study, contains
39~018
links, 12~982 nodes and 1~790 zones. We selected a set of potential road network
upgrades, listed in Table~\ref{tab:chicago_mods} and shown in
Figure~\ref{fig:chicago_mods}, based on data from the
Chicago Metropolitan Agency for Planning Transportation Improvement
Program website \citep{tip}.
In this set of upgrades, we will assume, based on the considerations
described in Section~\ref{sec:pairwise_predictor}, that only
07-94-0027 and 07-96-0013 will have significant pairwise interaction,
since they are two extensions to the same road.

\begin{table}
  \centering
  \small
  \begin{tabular}{llrl}
    \hline
    Project Id & Type & Cost (\$000s) & Description \\
    \hline
    03-02-9005 & capacity upgrade & 999 & add lanes on I-90 to I-294 \\
    03-03-0101 & capacity upgrade & 465 & add lanes on Meacham Rd.\\ 
    03-95-0001 & new road & 4000 & Elgin-O'Hare Expwy to Lake St. \\
    03-96-0024 & capacity upgrade & 1000 & widen lanes on Villa St./Lake St. \\
    07-06-0014 & new road & 472 & from Cottage Grove Ave. to Mark Collins Dr. \\
    07-94-0027 & new road & 700 & Joe Orr Rd. extension to Burnham Ave. \\
    07-96-0013 & new road & 748 & Joe Orr Rd. extension to Sheffield Ave. \\
    07-97-0055 & capacity upgrade & 4000 & add lanes on I-57 from I-80 to I-294 \\
    \hline
  \end{tabular}
  \caption{Potential upgrades to the Chicago Regional road network.}
  \label{tab:chicago_mods}
\end{table}

\begin{figure}
  \begin{center}
  \includegraphics[width=\textwidth]{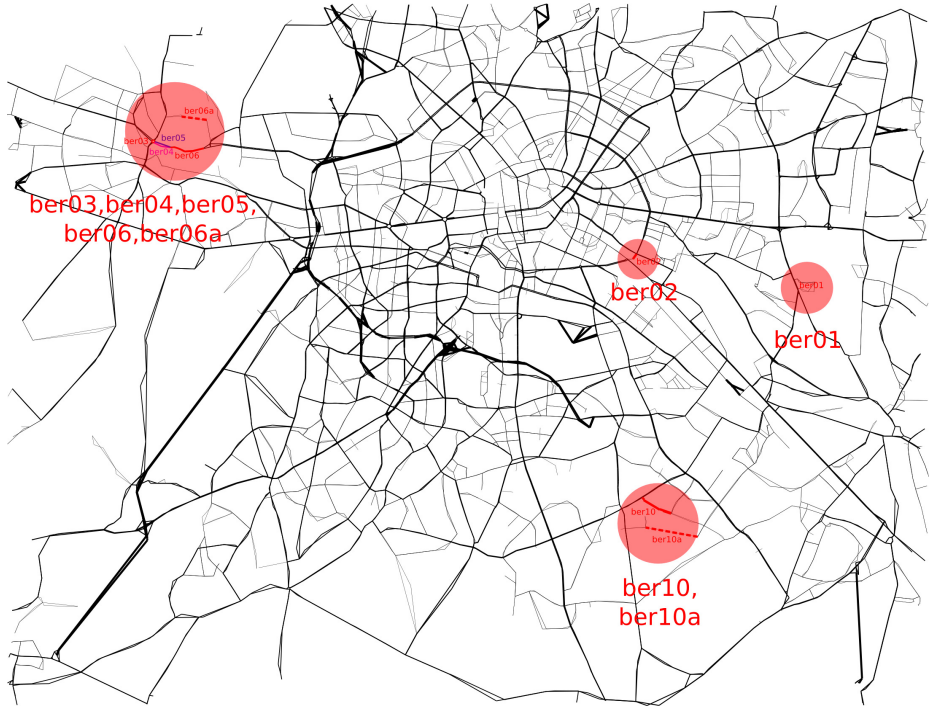}
  \end{center}
  \caption{Hypothetical upgrades to Berlin Center road network.}
  \label{fig:berlin_mods_overview}
\end{figure}

\begin{figure}
  \begin{center}
  \includegraphics[width=0.5\textwidth]{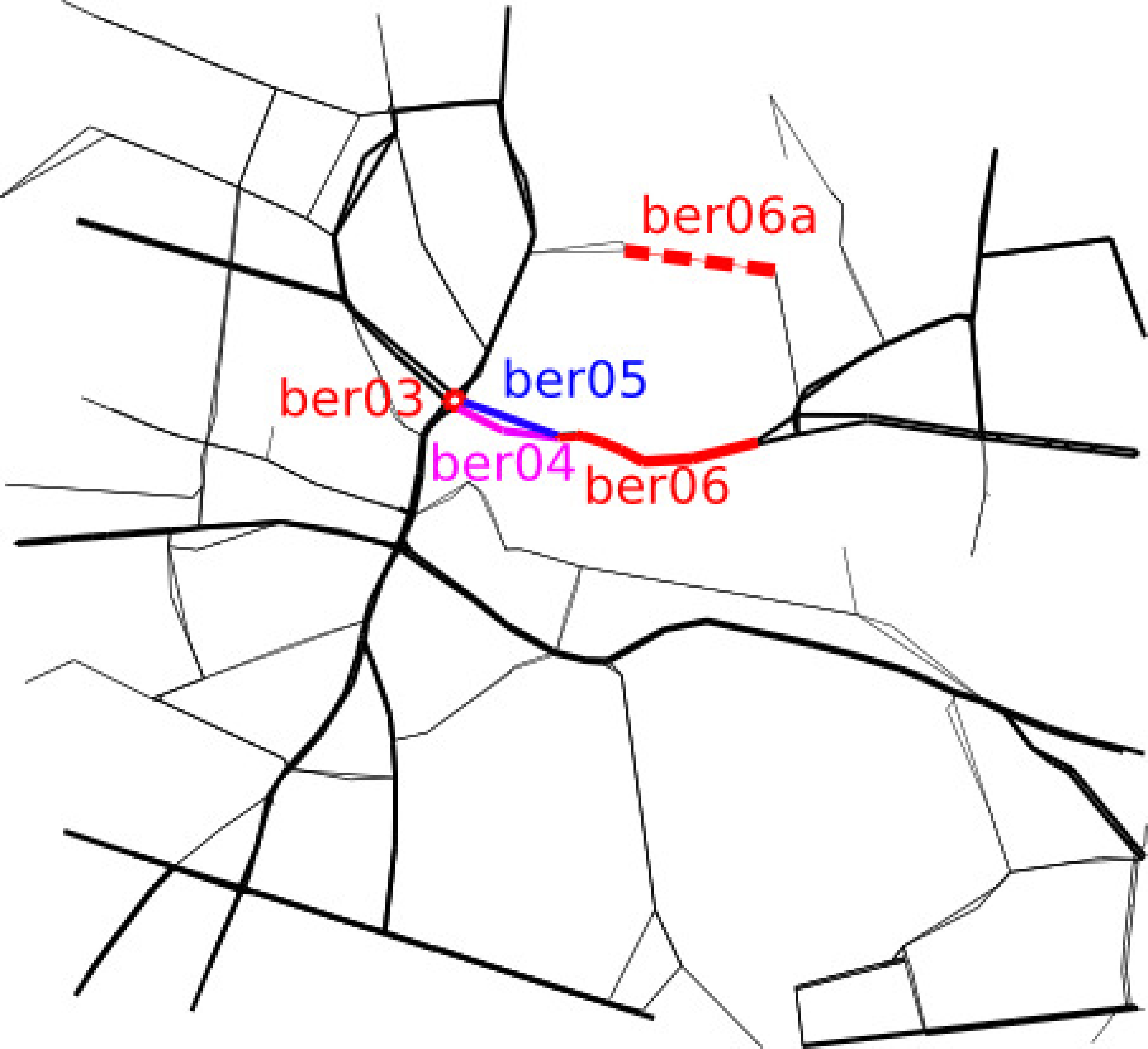}
  \end{center}
  \caption{Detail showing a set of five significantly interacting hypothetical upgrades in the Berlin Center road network.}
  \label{fig:berlin_mods_ber3456}
\end{figure}

\begin{figure}
  \begin{center}
  \includegraphics[width=0.3\textwidth]{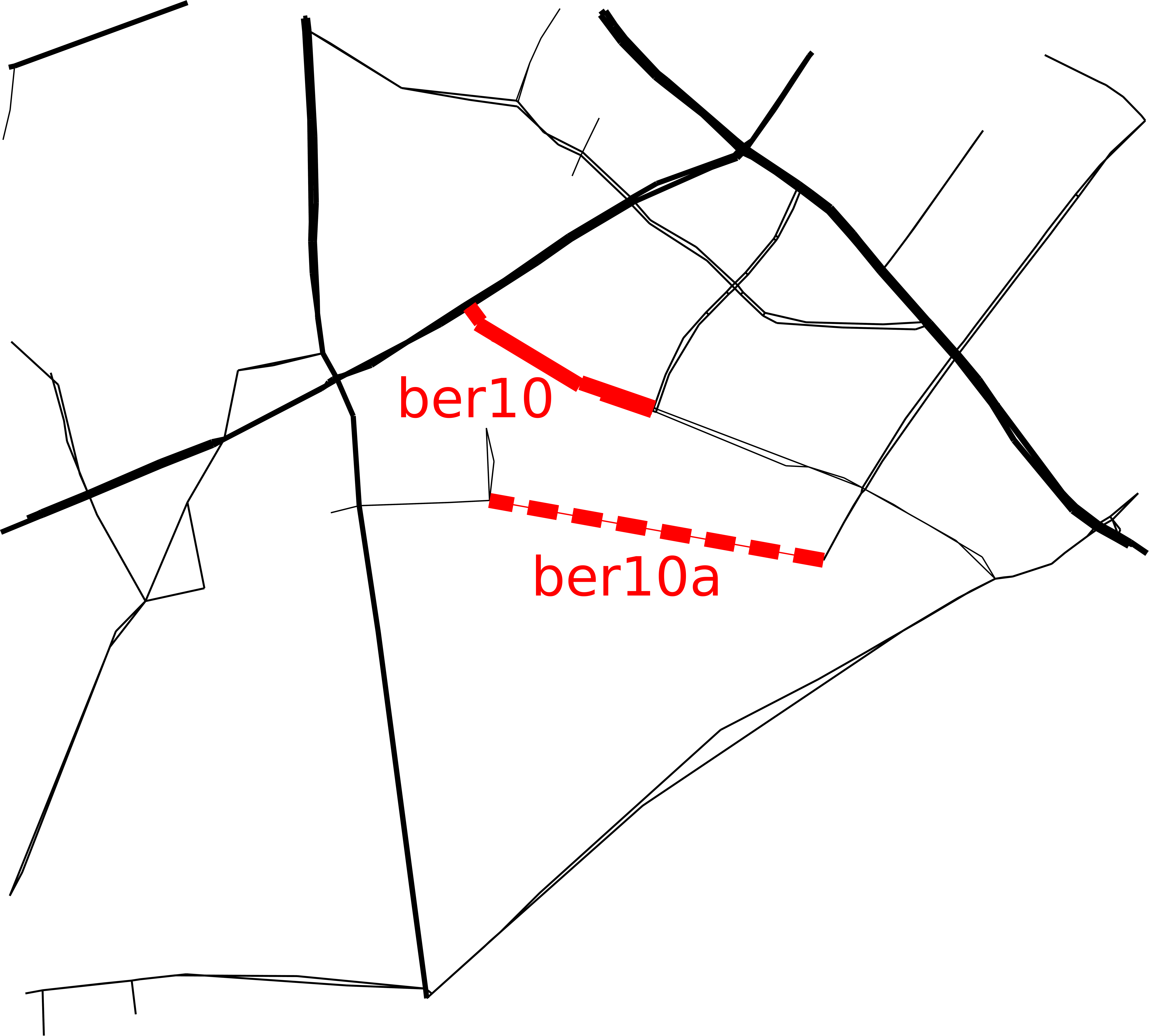}
  \end{center}
  \caption{Detail showing a set of two significantly interacting hypothetical upgrades in the Berlin Center road network.}
  \label{fig:berlin_mods_ber10}
\end{figure}

The second road network is the Berlin Center road network, used in
\citet{jahn05}. This network consists of 28~376 links, 12~981 nodes
and 865 zones. We created nine hypothetical road network upgrades in this
network by increasing capacity or adding roads in areas that show traffic
congestion when the TAP is solved on the network with all the O-D demands
from the original data doubled. These upgrades are listed in
Table~\ref{tab:berlin_mods} and shown in Figure~\ref{fig:berlin_mods_overview}.
Figures \ref{fig:berlin_mods_ber3456} and \ref{fig:berlin_mods_ber10}
are enlargements of the two clusters of nearby upgrades, in which 
each pair of upgrades is considered to be have significant pairwise
interaction according to the methods described in Section~\ref{sec:pairwise_predictor}.

\begin{table}
  \centering
  \begin{tabular}{llr}
    \hline
    Project Id & Type & Cost (\$000s) \\
    \hline
    ber01 & capacity upgrade & 300   \\
    ber02 & capacity upgrade & 1000  \\
    ber03 & capacity upgrade & 800   \\
    ber04 & capacity upgrade & 2500  \\
    ber05 & capacity upgrade & 2000  \\
    ber06 & capacity upgrade & 4000    \\
    ber06a & new road & 8000   \\
    ber10  & capacity upgrade &  1200  \\
    ber10a & new road & 8700   \\
    \hline
  \end{tabular}
  \caption{Hypothetical upgrades to the Berlin Center road network.}
  \label{tab:berlin_mods}
\end{table}

We used our implementation of the Frank-Wolfe method to solve the TAP
on the baseline (unmodified) Chicago Regional roads network, giving
the baseline total travel time of 33~657~132 VHT.
This was computed after 10~000
iterations, giving a relative gap of $1.2\times10^{-5}$. All
subsequent TAP runs were computed to the same relative gap.
For the Berlin Center road network, the baseline total travel time is
20~817~229 VHT, and all TAP runs are solved to a relative gap
of $1.0 \times 10^{-6}$.

The upgrade subset model (Equations
\ref{eqn:objective}--\ref{eqn:budget_constraint}) was written in
MiniZinc and solved with the \eclipse interval constraint solver.

In order to evaluate the accuracy of our method, we also computed
the VHT resulting from solving the TAP with our implementation 
for the road network modified with each of the 8 individual upgrades 
in Table~\ref{tab:chicago_mods}, as well as each of the 28 pairs
of upgrades (thereby obtaining every value of $d_{ij}$
as discussed in Section~\ref{sec:upgrade_subset_model}),
and each of the 219 subsets (with 3 or more elements) of the
set of potential upgrades. The latter 219 values
are used as the gold standard to evaluate the accuracy of estimating
the VHT change for a set of upgrades from only individual and pairwise
changes. This allows us to calculate a relative error
\begin{equation}
  \label{eqn:error}
  \frac{\left|\sum_{i \in S} v_i + \sum_{i \neq j \in S} d_{ij} - 
        \Delta\mbox{VHT}_S  \right|}
  {\Delta\mbox{VHT}_S}
\end{equation}
where $v_i$ and $d_{ij}$ are as defined in Section~\ref{sec:upgrade_subset_model}
giving the estimate of the VHT change used in the objective function
(Equation~\ref{eqn:objective}), and $\Delta\mbox{VHT}_S$ is the 
change in VHT computed by solving the TAP for the road network
modified by all upgrades in the
subset $S$ of the set of potential upgrades.
This can be extended to compute the relative error when accounting for
$k$-wise changes for all subsets up to (and including) size $k \leq \card{S}$.

\begin{table}
  \centering
\small
\begin{tabular}{lrrrr}
\hline
Data used &  $\Delta\mbox{VHT}$ computations  & Mean Error \% & Num Error $> 10\%$ \\
\hline
individual only  &  8  &  4.0  & 32  \\
 significant pairwise  &  9  &  2.4  & 3  \\
 all pairwise  &  36  &  0.68  & 2  \\
 all subsets size $\leq$ 3   &  92  &  0.25  & 1  \\
 all subsets size $\leq$ 4   &  162  &  0.12  & 0  \\
 all subsets size $\leq$ 5   &  218  &  0.043  & 0  \\
 all subsets size $\leq$ 6   &  246  &  0.0089  & 0  \\
 all subsets size $\leq$ 7   &  254  &  0.0015  & 0  \\
 all subsets size $\leq$ 8   &  255  &  0.0  & 0  \\
 \hline
\end{tabular}
  \caption{Relative error in Chicago Regional $\Delta\mbox{VHT}$ values 
    with different degrees of upgrade interactions (sizes of subsets
    of upgrades) computed.}
  \label{tab:chicago_error_stats}
\end{table}

Table~\ref{tab:chicago_error_stats} shows the relative errors in 
the Chicago Regional road network for the
set of upgrades in Table~\ref{tab:chicago_mods}.
Using all pairwise $\Delta\mbox{VHT}$ computations, we find that
the mean value of this relative error is 0.68\%. 

If we use only the individual VHT computations, \ie assuming all the
$d_{ij}$ are 0, thereby requiring only 8 TAP computations, rather than
36 for all pairwise interactions as before, then the mean value of this
relative error is 4.0\%, and 32 subsets have a relative error larger
than 10\% (compared to only 2 when using all pairwise interactions).

If, however, we add to the individual VHT computations only the single
pairwise interaction between 07-94-0027 and 07-96-0013, thereby
requiring only 9 TAP computations, the mean relative error is 2.4\%,
and only 3 subsets have a relative error greater than 10\%.

In summary, this small example shows that requiring only a single extra
TAP computation, that for a pair of upgrades that \latin{a priori} 
would be expected to have a significant difference between the sum
of their individual VHT improvements and the VHT improvement resulting from 
implementing both, reduces the relative error from 4\% to 2.4\%.
If we assume that the VHT improvement due to 
all pairwise interactions are significantly
different from the sum of their individual VHT improvements, the error is reduced
further, to 0.68\%,
however this 
requires a total of 36 TAP computations (an increase of 28 computations),
rather than just a single
extra computation.

\begin{table}
  \centering
\small
\begin{tabular}{lrrrr}
\hline
Data used &  $\Delta\mbox{VHT}$ computations  & Mean Error \% & Num Error $> 10\%$ \\
\hline
individual only  &  9  &  8.4  & 197  \\
 significant pairwise  &  20  &  1.8  & 16  \\
 all pairwise  &  45  &  1.8  & 16  \\
 all subsets size $\leq$ 3   &  129  &  0.24  & 0  \\
 all subsets size $\leq$ 4   &  255  &  0.038  & 0  \\
 all subsets size $\leq$ 5   &  381  &  0.019  & 0  \\
 all subsets size $\leq$ 6   &  465  &  0.0062  & 0  \\
 all subsets size $\leq$ 7   &  501  &  0.0013  & 0  \\
 all subsets size $\leq$ 8   &  510  &  0.00012  & 0  \\
 all subsets size $\leq$ 9   &  511  &  0.0  & 0  \\
 \hline
\end{tabular}
  \caption{Relative error in Berlin Center $\Delta\mbox{VHT}$ values 
    with different degrees of upgrade interactions (sizes of subsets
    of upgrades) computed.}
   \label{tab:berlin_error_stats}
\end{table}

We also used the 9 individual upgrades to the Berlin Center
network in Table~\ref{tab:berlin_mods}, along with all 36 pairs
of upgrades and 466 subsets of size 3 or greater to assess accuracy
in the Berlin Center network.
This example is more complex than the Chicago Regional example,
in that the set of five upgrades shown in Figure~\ref{fig:berlin_mods_ber3456}
have a higher degree of interaction than just pairwise.

The results for the Berlin Center road network are shown in 
Table~\ref{tab:berlin_error_stats}. Here we can see that the mean
error of 8.4\% when only individual upgrade TAP computations are used
is reduced to 1.8\% when all pairwise computations are used, requiring
45 TAP computations. However if only an additional 11 TAP computations
(for a total of 20) are made for the two sets of significantly interacting
upgrades (shown in Figure~\ref{fig:berlin_mods_ber3456} and 
Figure~\ref{fig:berlin_mods_ber10}) then the mean error is reduced to 1.8\%,
the same as using all pairwise computations.
This example is more complex than the Chicago Regional example,
in that the set of five upgrades shown in Figure~\ref{fig:berlin_mods_ber3456}
have a higher degree of interaction than just pairwise. As a consequence,
if we use only all subsets of the interacting upgrades (those in Figures
\ref{fig:berlin_mods_ber3456} and \ref{fig:berlin_mods_ber10}), then
the mean error is reduced to 0.02\% with 129 $\Delta\mbox{VHT}$ computations.

Table~\ref{tab:chicago_error_stats} and Table~\ref{tab:berlin_error_stats}
also show that we can obtain more accuracy by solving the TAP for all subsets
of increasing size, up to subset size equal 
to the number of individual upgrades, at which point the TAP is solved
for all subsets and the error is necessarily zero. However going beyond
pairwise interaction quickly leads to an impractical number of TAPs
to solve, for diminishing returns. In addition, if we 
used higher degrees of interaction than pairwise then the
objective function (Equation~\ref{eqn:objective}) would
accordingly have a degree
higher than quadratic.

\begin{table}
  \centering
\small
\begin{tabular}{lrrrr}
\hline
Data used &  $\Delta\mbox{VHT}$ computations  & Mean Error \% & Num Error $> 10\%$ \\
\hline
individual only  &  100  &  0.66  & 1  \\
 significant pairwise  &  145  &  0.22  & 0  \\
 all pairwise  &  5050  &  0.0  & 0  \\
 \hline
\end{tabular}
   \caption{Error relative to using all pairwise upgrades
            in Berlin Center $\Delta\mbox{VHT}$ values when using only
            individual upgrades versus only significant upgrade pairs.}
   \label{tab:berlin_100_pairwise_error_stats}
\end{table}

In order to show the effectiveness of this approach on a larger scale,
we generated 100 random upgrades (20 new road links
and 80 roads with increased capacity) in the Berlin Center road network.
Of these 100 upgrades (and therefore 4~950 pairs of upgrades), 45 pairs
are considered to have potentially significant interaction by using
the same distance threshold used for the Berlin Center upgrades
from Table~\ref{tab:berlin_mods}.
We measured the error in $\Delta\mbox{VHT}$ value for each of the
pairs of upgrades when solving the TAP for individual upgrades only
versus solving the TAP additionally for only the upgrade pairs
considered significant.  
The results are shown in Table~\ref{tab:berlin_100_pairwise_error_stats},
where it can be seen that including the 45 potentially significant pairwise
computations reduces the mean error by a factor of three,
and results in no pairs with a relative error  greater than 10\%.
Note that the relative error used here is different from those used previously:
due to the infeasibly large number of subsets involved, the error is measured
relative not to every possible subset, but only to every pair (subset of size 2).

If we solve the model described in Section~\ref{sec:upgrade_subset_model}
with the Chicago Regional road network upgrades from 
Table~\ref{tab:chicago_mods} and a budget of \$10~000~000, we
find that the optimal subset of upgrades is
03-95-0001, 07-06-0014, 07-94-0027, 07-96-0013, and 07-97-0055,
with net value of 168 358 (thousand dollars).
There is a VHT reduction of 48~844 according to the model (using only the 8
individual and additional one significant pairwise $\Delta\mbox{VHT}$
computations), and the actual VHT reduction for this subset of upgrades
is 49~794, so the model's value has a relative error of 1.9\%.
Solving the model including all 3-wise interactions, thereby requiring
92 VHT computations and a cubic objective function, gives a VHT 
reduction of 49~927, a relative error of 0.27\%.

If we solve the model
with the hypothetical upgrades to the Berlin Center road network from
Table~\ref{tab:berlin_mods} and a budget of \$10~000~000, we
find the optimal subset of upgrades is ber01, ber06a, and ber10, 
with net value  456~532  (thousand dollars).
There is a VHT reduction of 127~679 according to the model (which is using
only the 9 individual and additional 11 significant pairwise 
$\Delta\mbox{VHT}$ computations). The actual VHT reduction computed for this
subset of upgrades is 127~653, so the model's value has an error of only
0.02\%. Solving the model using all 3-wise interactions, which requires
129 VHT computations, is enough to obtain a zero relative error in the VHT reduction for this optimal subset.

\section{Finding an optimal schedule of upgrades}
\label{sec:upgrade_schedule}

So far we have only considered the problem of finding the optimal
subset of a set of potential upgrades, to fit in a given budget.
However, a more realistic (and more difficult) problem is to find an
optimal schedule of upgrades over a period of time, in which the
upgrades implemented in each time period must fit within a budget for
that time period. We will define ``optimal'' in this context shortly,
but an example of a feasible schedule should help clarify the basic
idea. Given the set of potential upgrades from
Table~\ref{tab:chicago_mods}, an example of a feasible (but not 
necessarily optimal) schedule is
shown in Table~\ref{tab:schedule_example}.  In each time period (which
may be, for instance, a year, or a 5 year period), a set of upgrades
that fit within that period's budget is scheduled to be built.

\begin{table}
  \centering
  \begin{tabular}{rlrrrrrr}
  \hline
  \multicolumn{2}{c}{Time period $t$} & 1 & 2 & 3 & 4 & 5 & Total  \\
  \multicolumn{2}{c}{Budget (\$000s)}  & 1000 & 4000 & 1500 & 3000 & 5000 & 14500\\
  $i$ & Project Id  \\
  \hline
  1 & 03-02-9005 & X    &         &         &         &      \\
  2 & 03-03-0101 &      &         &         &    X    &      \\
  3 & 03-95-0001 &      &    X    &         &         &      \\
  4 & 03-96-0024 &      &         &         &    X    &      \\
  5 & 07-06-0014 &      &         &         &    X    &      \\
  6 & 07-94-0027 &      &         &    X    &         &      \\
  7 & 07-96-0013 &      &         &    X    &         &      \\
  8 & 07-97-0055 &      &         &         &         &    X \\
  \hline
  \multicolumn{2}{c}{Expenditure (\$000s) } & 999 & 4000 & 1448 & 1937 & 4000 & 12384 \\
  \hline
  \end{tabular}
  \caption{Example of a feasible road network upgrade schedule.}
  \label{tab:schedule_example}
\end{table}

We will define optimality in this context as the schedule that gives
the maximum \firstuse{net present value} of VHT reductions. Recall that
the present value of a monetary amount with value $c$ at time $t$ 
in the future is
$\frac{c}{(1+r)^t}$ where $r$ is the interest rate. We will assume here,
conventionally,
that our time periods are years and the interest rate $r$ is annual, so $t$
is in years. 

Using the same notation as previously (Section~\ref{sec:upgrade_subset_model}),
but with the addition of another subscript $t$ on variables to denote
their value at time $t \in \{1, 2,  \ldots T\}$ where $T$ is the number of time
periods, the optimal schedule can be determined by
the following constrained optimization problem:
\begin{equation}
\label{eqn:schedule_objective}
\max \left( 
\sum_{1 \leq t \leq T}
\sum_{1 \leq i < j \leq N} \frac{y_{it} y_{jt}  d_{ijt} m}{(1+r)^t} + 
\sum_{1 \leq i \leq N} y_{it} \left( \frac{v_{it} m}{(1+r)^t} - c_i \right)
\right)
\end{equation}
subject to
\begin{eqnarray}
  \label{eqn:schedule_budget_constraint}
  \sum_{1 \leq i \leq N} y_{it} c_i \leq B_t, && 1 \leq t \leq T
\end{eqnarray}
\begin{eqnarray}
  \label{eqn:schedule_onceonly_constraint}
  \sum_{1 \leq t \leq T} y_{it} \leq 1, && 1 \leq i \leq N
\end{eqnarray}

The first constraint (Equation~\ref{eqn:schedule_budget_constraint})
ensures that the total
expenditure in each year is within that year's budget. For convenience
we can assume that all budgets and costs are expressed as their present value.
The second constraint (Equation~\ref{eqn:schedule_onceonly_constraint})
states that each
potential upgrade is built at most once. The example of a feasible schedule
in Table~\ref{tab:schedule_example} illustrates these constraints:
the expenditure (sum of project costs) in each column is less than or equal to the budget for that
time period, and each row has only one cross indicating that the
project is built in only one time period.

Solving this problem exactly therefore requires not just all the
individual VHT changes $v_i$ and pairwise changes $d_{ij}$ as in
Section~\ref{sec:upgrade_subset_model}, but these values at all time
periods $1 \leq t \leq T$. They may have different values
depending on which upgrades have been already built at
earlier time periods, due to the effects of upgrades.  In addition,
although we have assumed that the O-D demand matrix is constant in a
given time period, presumably it will change over time, for example
increased demand generally and in between particular origins and
destinations (new suburbs, for example). Hence finding the exact
solution of this problem is impractical due to the factorial number of
TAPs that would need to be solved (potentially we would need to solve
the TAP for every permutation of feasible upgrade subsets).

We therefore propose the following heuristic ``greedy'' strategy
to find a feasible schedule that attempts to approximate the optimal
schedule.  First, we choose the subset of upgrades that is optimal
at the end of the entire planning period, within the total budget over all 
time periods. That is,
\begin{equation}
\label{eqn:final_pv_subset_objective}
\max \left( 
\sum_{1 \leq i < j \leq N} \frac{y_{i} y_{j} d_{ijT} m}{(1+r)^T} + 
\sum_{1 \leq i \leq N} y_{i} \left( \frac{v_{iT} m}{(1+r)^T} -c_i \right)
\right)
\end{equation}
subject to
\begin{equation}
  \label{eqn:final_pv_budget_constraint}
  \sum_{1 \leq i \leq N} y_i c_i \leq \sum_{1 \leq t \leq T} B_t
\end{equation}
giving a set of upgrades $U$ that will all be built by the end of the
planning period. Then we need to schedule the building of these
upgrades in the planning period to maximize their net present value
(in terms of VHT reduction). To do so,
at each time period $t \in 1,2, \ldots .. T$ solve the optimal
upgrade subset problem (Section~\ref{sec:upgrade_subset_model},
Equations~\ref{eqn:objective}--\ref{eqn:budget_constraint}), modified
however to use the net present value as is done in Equation~\ref{eqn:schedule_objective}. That is, instead of the objective function Equation~\ref{eqn:objective},
use:

\begin{equation}
\label{eqn:npv_subset_objective}
\max \left( 
\sum_{i, j\in U_t, i < j} \frac{y_{it} y_{jt} d_{ijt} m}{(1+r)^t} + 
\sum_{i \in U_t} y_{it}  \left( \frac{ v_{it} m}{(1+r)^t } - c_i \right)
\right)
\end{equation}
where $U_t$ is the current set of potential upgrades (at $t=1$, $U_t$ is the
optimal subset of upgrades $U$ for the end of the planning
period obtained from 
solving the constrained optimization problem given by
Equations \ref{eqn:final_pv_subset_objective}
and \ref{eqn:final_pv_budget_constraint}), and the budget $B$ 
in Equation~\ref{eqn:budget_constraint} is $B_t$ \ie the budget for time period $t$ is used. 

Then fix the values of $y_{ijt}$ for the current $t$, 
giving the subset of upgrades
for time period $t$ and a new (upgraded) road network, and a new
set of potential upgrades $U_{t+1} = U_t \setminus V$ where $V \subseteq U_t$
is the set of upgrades just chosen (\ie those $i$ for which $y_{ijt} =1$).
So at each  time period $1 < t \leq T$,
the optimal upgrade subset problem is solved for the set of upgrades
not yet built, using the road network determined by
the upgrades chosen in period $t -1$ (and transitively by all previous
time periods), and the O-D demand matrix for time $t$.

In terms of the example shown in Table~\ref{tab:schedule_example},
this heuristic algorithm can be intuitively thought of as simply
choosing the optimal subset (of the set of upgrades found as optimal
at the final time period given the total budget)
in each time period (column) starting from
the first (leftmost column), assuming in each column that the road
network now has all upgrades made in previous time periods, and hence
they cannot be chosen again, and that these choices are irrevocable.

When using this heuristic algorithm, we use the same values of
$d_{ij}$ as in Section~\ref{sec:upgrade_schedule}, for each time
period.  That is, we assume $d_{ij} = 0$ unless, for the reasons
discussed in Section~\ref{sec:pairwise_predictor}, there is thought to
be a significant difference between the sum of the VHT changes for
the individual upgrades $i$ and $j$ and their pairwise VHT change,
in which case $d_{ij}$ is computed by solving the TAP under the
conditions prevailing at the current time period (\ie the network
including upgrades at previous time periods, and with the current O-D
demand matrix).

As an alternative to the greedy heuristic, we could make stronger assumptions
about the independence of the individual upgrades
to reduce the complexity of the optimal upgrade schedule problem.
For example, if all the upgrades are assumed to not have significant
interaction (\ie all the $d_{ij}$ are assumed to be 0), then the
quadratic term in Equation~\ref{eqn:schedule_objective} disappears
and the TAP only needs to be solved for each individual upgrade
at each time period. However, we would expect, given the results
discussed in Section~\ref{sec:results}, that ignoring these 
interactions entirely could lead to significant errors in VHT reduction
estimates and therefore potentially very far from optimal schedules.

These algorithms demonstrate a shortcoming of the Frank-Wolfe method
that we have been using. Some algorithms such as Origin-Based
Assignment (OBA) \citep{bargera02} and Algorithm B \citep{dial06} can
``warm start'' from a converged TAP solution with a changed O-D demand
matrix and reach equilibrium more quickly than if started from
scratch. Frank-Wolfe does not have this ability and hence OBA or
Algorithm B may be useful in the situation where we solve the TAP for
a number of modified networks with a given O-D demand matrix, and then
need to re-solve it with a changed O-D demand matrix in order to find the baseline
VHT for the next scheduling period. They will not necessarily be an 
improvement in the majority of cases where we need to solve the TAP
again with a modified network, but obviously any improvement in the speed
with which the TAP can be solved is useful.

\subsection{Results}
\label{sec:schedule_results}

We implemented the heuristic algorithm by writing the models
(Equations \ref{eqn:final_pv_subset_objective}--\ref{eqn:npv_subset_objective})
in MiniZinc and solving with the \eclipse interval constraint solver,
and using a scripting language to run the solvers with the appropriate
data at each iteration.

We used the eight upgrades to the Chicago Regional
road network listed in Table~\ref{tab:chicago_mods}, with the five time
periods and budgets for each as in the example schedule
in Table~\ref{tab:schedule_example}.
We used the Chicago Regional data, identical to that previously used
in the first time period. For the subsequent time periods, we 
re-ran the TAP for all modifications with a 10\% increase in trips
to and from the CBD and a 15\% increase in trips to and from
O'Hare airport in each period.  We used the factor $m = 3650$
to convert VHT to dollars and an annual interest rate of 
4\%.
The resulting schedule is shown in 
Table~\ref{tab:chicago_heuristic_schedule}. 
If all the upgrades
were built at the end of the planning period (\ie the first
step of the greedy heuristic, where we find a subset that is
optimal at the end, within the total budget over all periods),
then the dollar present value of the VHT reduction (in thousands of dollars)
is 148~151,
while the net present value of the schedule shown in 
Table~\ref{tab:chicago_heuristic_schedule} is 164~359.

\begin{table}
  \centering
\begin{tabular}{rlrrrrrr}
\hline
\multicolumn{2}{c}{Time period $t$} & 1 & 2 & 3 & 4 & 5 & Total  \\
\multicolumn{2}{c}{Budget (\$000s)}  & 1000 & 4000 & 1500 & 3000 & 5000 & 14500\\
$i$ & Project Id  \\
\hline
1 & 03-02-9005 &      &         &         &         &      \\
2 & 03-03-0101 &      &         &         &    X    &      \\
3 & 03-95-0001 &      &    X    &         &         &      \\
4 & 03-96-0024 &      &         &         &         &      \\
5 & 07-06-0014 &      &         &    X    &         &      \\
6 & 07-94-0027 &      &         &    X    &         &      \\
7 & 07-96-0013 & X    &         &         &         &      \\
8 & 07-97-0055 &      &         &         &         &    X \\
\hline
\multicolumn{2}{c}{Expenditure (\$000s) } & 748 & 4000 & 1172 & 465 & 4000 & 10385 \\
\hline
\end{tabular}
  \caption{Road network upgrade schedule for the Chicago Regional road
     network generated by the greedy heuristic algorithm.}
   \label{tab:chicago_heuristic_schedule}
\end{table}

For comparison, we  used the same data in the simplified model that 
assumes independence of all upgrades. As previously mentioned,
this means that the quadratic term in Equation~\ref{eqn:schedule_objective}
disappears, giving the following linear constrained optimization problem:
\begin{equation}
\label{eqn:npv_independent_schedule_objective}
\max \left( 
\sum_{1 \leq t \leq T}
\sum_{1 \leq i \leq N} y_{it} \left( \frac{v_{it} m}{(1+r)^t} - c_i \right)
\right)
\end{equation}
subject to
\begin{eqnarray}
  \label{eqn:independent_schedule_budget_constraint}
  \sum_{1 \leq i \leq N} y_{it} c_i \leq B_t, && 1 \leq t \leq T
\end{eqnarray}
\begin{eqnarray}
  \label{eqn:independent_schedule_onceonly_constraint}
  \sum_{1 \leq t \leq T} y_{it} \leq 1, && 1 \leq i \leq N
\end{eqnarray}
It is practical to solve this exactly (again by writing
the model in MinZinc and solving with the \eclipse interval constraint solver)
as the assumption of the independence of all the upgrades means we need
only solve the TAP for each individual upgrade at each time period
(each subsequent time period having increased O-D demands as described
above). The resulting upgrade schedule is shown 
in Table~\ref{tab:npv_chicago_independent_schedule}. 

This schedule has a net present value (in thousands of dollars)
of 145~648 according to the model, and actual NPV (computed by solving
the TAP for each upgrade at each period according to the schedule and
then computing the NPV of the upgrades)
of 164~091.

\begin{table}
  \centering
\begin{tabular}{rlrrrrrr}
\hline
\multicolumn{2}{c}{Time period $t$} & 1 & 2 & 3 & 4 & 5 & Total  \\
\multicolumn{2}{c}{Budget (\$000s)}  & 1000 & 4000 & 1500 & 3000 & 5000 & 14500\\
$i$ & Project Id  \\
\hline
1 & 03-02-9005 &      &         &         &         &      \\
2 & 03-03-0101 &      &         &    X    &         &      \\
3 & 03-95-0001 &      &    X    &         &         &      \\
4 & 03-96-0024 &      &         &         &         &      \\
5 & 07-06-0014 &      &         &         &    X    &      \\
6 & 07-94-0027 &      &         &         &         &    X \\
7 & 07-96-0013 & X    &         &         &         &      \\
8 & 07-97-0055 &      &         &         &         &    X \\
\hline
\multicolumn{2}{c}{Expenditure (\$000s) } & 748 & 4000 & 465 & 472 & 4700 & 10385 \\
\hline
\end{tabular}
  \caption{Road network upgrade schedule for the Chicago Regional road
     network generated by solving the linear model assuming independence
     of all upgrades.}
   \label{tab:npv_chicago_independent_schedule}
\end{table}

Both of these schedules are improvements on the example schedule
in Table~\ref{tab:schedule_example}, however, which has a a net present
value of 157~676 (thousand dollars).

  
This example demonstrates that the greedy heuristic algorithm, using
information about pairwise interaction between upgrades, is
capable of finding better upgrade schedules than the exact solution
of the simplified upgrade scheduling problem in which pairwise interactions
are ignored.

\section{Conclusions}

In this paper, we have demonstrated a method to find the optimal subset
of a set of potential road network upgrades, without having to 
solve the traffic assignment problem for each subset. Although
this may be necessary in principle, we have shown that reasonable
approximations can be achieved by considering only those pairs of
upgrades that potentially have significant interactions. This subset of
pairs  can be predicted with a simple distance threshold or clustering
technique. In addition, we have shown how to find an implementation
schedule of the upgrades to maximize the net present value of the reduction
in vehicle hours travelled over the planning period.

It should be noted that many major simplifications have been made in
this modelling process. Perhaps most importantly, we have assumed that
the single optimization objective is to minimize vehicle hours
travelled, a value determined entirely by the flow on the road
network at user equilibrium, given a fixed origin-destination demand
matrix.  However there may be important other benefits, or drawbacks,
to a particular road upgrade, such as social and environmental costs
(or benefits), or financial or political considerations relating to
building roads or operating toll roads for example, that are not
accounted for. If, as is done in Section~\ref{sec:upgrade_schedule},
 the VHT value is converted to a monetary value (by
assigning a monetary cost to each vehicle hour --- and here different
types of traffic such as personal, work-related, commercial freight,
etc. should be distinguished), and if externalities are also
converted to monetary costs, then they could be incorporated as costs
or benefits in a more sophisticated model, however that is beyond the
scope of this work.

We have also simplified somewhat by considering only road traffic
assignment, ignoring mode choices where such alternatives as public
transport exist.  Traffic assignment is usually used as one step in an
iterative multi-step process including mode choice and trip generation:
this multi-step process could be used in place of the simple traffic
assignment described here to incorporate different transport modes 
into the model.

By using a conventional static traffic assignment model, where flow on
a link is allowed to exceed its capacity when congested, we have
ignored the possibility of demand being affected by congestion. Some
more sophisticated models can address this, such as ``stationary
dynamic'' models \citep{nesterov03} in which excess demand is queued,
or elastic demand models \citep{babonneau08}, in which demand for
travel between an origin and destination is a function of travel
costs.  Further, by using a static model with fixed origin-destination
demand, we have ignored any induced demand that may result from
improving travel times, or at least deferred this to the separate
problem of estimation of demand in the future.

\section{Acknowledgements}

We have been greatly assisted in this work by the public availability
of the transportation network test problems, TAP file format, and
sample TAP code from Dr Hillel Bar-Gera. This research made use of the
Victorian Partnership for Advanced Computing HPC facility and support
services. NICTA is funded by the Australian Government as represented
by the Department of Broadband, Communications and the Digital Economy
and the Australian Research Council through the ICT Centre of
Excellence Program.


\end{document}